\documentclass[11pt,twoside, reqno]{amsart}

\usepackage{amsmath}
\usepackage{amsthm}
\usepackage{amsfonts, amssymb}
\usepackage{mathrsfs}
\usepackage[all]{xy}
\usepackage{url}
\usepackage{tikz}

\usepackage{latexsym}

\theoremstyle{plain}
\newtheorem{lema}{Lemma}
\newtheorem{prop}[lema]{Proposition}
\newtheorem{teo}[lema]{Theorem}

\newtheorem*{intro1}{Corollary \ref{maincoro}}
\newtheorem*{intro2}{Corollary \ref{maincoro2}}
\newtheorem*{teomain}{Theorem \ref{main}}
\newtheorem{coro}[lema]{Corollary}
\theoremstyle{remark}

\newtheorem{obs}[lema]{Remark}

\theoremstyle{definition}
\newtheorem{defi}[lema]{Definition}
\newtheorem{ej}[lema]{Example}

\newcommand{\x}{\mathcal{X}}
\newcommand{\kp}{\mathcal{K}}

\newcommand{\N}{\mathbb{N}}
\newcommand{\R}{\mathbb{R}}
\newcommand{\Z}{\mathbb{Z}}

\newcommand\restr{\raisebox{-0.3ex}{$|$}\raisebox{0.3ex}{}}
\newcommand{\aatop}{\mathsf{ATop}}
\newcommand{\tatop}{\mathsf{T_0ATop}}
\newcommand{\topp}{\mathsf{Top}}

\usepackage{float}
\usepackage{graphicx}
\usepackage{xfrac}
\usepackage{adjustbox}
\usepackage{amsmath}
\usepackage{amsthm}
\usepackage{quiver}
\usepackage{soul}
\usepackage{graphics}
\usepackage{tikz}
\usepackage{tikz-cd}
\usetikzlibrary{babel}
\usepackage[english]{babel}
\usepackage{hyperref}
\hypersetup{
	colorlinks=true,
	linkcolor=blue,
	filecolor=magenta,      
	urlcolor=cyan,
	pdftitle={Overleaf Example},
	pdfpagemode=FullScreen,
}
\usepackage{color}

\newcommand{\C}{\mathbb{C}}

\theoremstyle{remark}

\theoremstyle{definition}

\DeclareEmphSequence{\color{blue}\itshape}
\begin{document}
	
\pagestyle{myheadings} \markboth{{\sc  j.a. barmak and j. galotto}}{{\sc Directed systems of A-spaces and weak cores}}

\title[Directed systems of A-spaces and weak cores]{Directed systems of A-spaces and weak cores}

\author[J.A. Barmak]{Jonathan Ariel Barmak}
\author[J. Galotto]{Juan Galotto}

\thanks{Both authors are partially supported by CONICET and by grant UBACYT 20020250100097BA}


\address{Universidad de Buenos Aires. Facultad de Ciencias Exactas y Naturales. Departamento de Matem\'atica. Buenos Aires, Argentina.}

\address{CONICET-Universidad de Buenos Aires. Instituto de Investigaciones Matem\'aticas Luis A. Santal\'o (IMAS). Buenos Aires, Argentina. }

\email{jbarmak@dm.uba.ar}

\email{juangalotto@gmail.com}

\begin{abstract}
We prove that for any regular CW-complex $X$ there exists a quotient map $X\to Y$ that is a weak homotopy equivalence and which is minimal with that property, meaning that any quotient $Y\to Z$ that is a weak equivalence, must be a homeomorphism. The space $Y$ we obtain is an Alexandroff space: a space in which arbitrary intersections of open sets are open. The weak core $Y$ of $X$ is not unique in general. We study homotopy invariants preserved by filtered colimits of Alexandroff spaces. This can be seen as an analogous to Stong's theory of cores of finite topological spaces.

\end{abstract}

\makeatletter
\@namedef{subjclassname@2020}{%
  \textup{2020} Mathematics Subject Classification}
\makeatother



\maketitle

\section{Introduction} \label{intro}

An Alexandroff space (or A-space) is a topological space in which arbitrary intersections of open sets are open. In 1937, Alexandroff \cite{Ale} showed that A-spaces are the same as preordered sets. More specifically, there is an isomorphism between the category $\aatop$ of A-spaces with continuous functions, and the category of preordered sets with order-preserving maps, which is the identity on the underlying sets and the corresponding maps. Moreover, this functor restricts to an isomorphism between $T_0$ A-spaces and posets.

Every topological space is weak homotopy equivalent to an A-space by a result by McCord \cite{Mcc}: for every simplicial complex $K$ there exists an A-space $\mathcal{X}(K)$ and a weak equivalence $|K|\xrightarrow{\mu}\mathcal{X}(K)$. Also, for every A-space $X$ there is a simplicial complex $\mathcal{K}(X)$ and a weak equivalence $|\mathcal{K}(X)|\xrightarrow{\mu}X$. Moreover, $\mathcal{X}$ and $\mathcal{K}$ are in fact functors between $\aatop$ and the category of simplicial complexes. One of the interesting features of A-spaces is that a weak equivalence between A-spaces need not be a homotopy equivalence. Homotopy equivalences in $\aatop$ correspond to simplicial morphisms with an inverse modulo contiguity \cite{BM}.

An important class of A-spaces is given by finite topological spaces. Homotopies and homotopy types in finite spaces are elegantly described by a result of Stong \cite{Sto}. The core of a finite space $X$ is a strong deformation retract which has no proper strong deformation retracts. Equivalently, the core of $X$ is a space $Y$ such that there is a quotient map $X\to Y$ which is also a homotopy equivalence and where $Y$ is minimal with that property. Every finite space has a core and it is unique up to homeomorphism. 
The core $Y$ of $X$ can be obtained by a sequence of retractions $X=X_0\to X_1\to \ldots \to X_n=Y$, each of which maps any point $x$ to a point comparable with $x$. This idea has been extended to A-spaces in the following way. Suppose $X$ is an A-space and $X=X_0\xrightarrow{r_0} X_1\xrightarrow{r_1} \ldots $ is a sequence of comparative retractions, in the sense above. If every trajectory $x,r_0(x),r_1r_0(x),\ldots$ stabilizes, then $X$ is homotopy equivalent to the subspace of stable points \cite{Kuk}. One instance of this result is the Li-Milner Structure Theorem \cite{LM} which guarantees the existence of a core for A-spaces $X$ satisfying certain conditions. Moreover, under these conditions, the core is unique up to homeomorphism. A core here is understood again as the codomain of a homotopy equivalence that is a quotient, and which is minimal with this property. The aim of these generalizations is to find a minimal quotient that retains some desired invariants of our original space (e.g. the fixed point property). When not every trajectory stabilizes, the set of stable points gives little information. In fact, this set could be empty. In \cite{Gal} we propose then to consider the space of trajectories, which is a quotient of the original space $X$. It turns out that this space is not homotopy equivalent to $X$ in general, but weak homotopy equivalent. This idea resembles that of Morse theory, in which the gradient flow of a Morse function codifies deformations which induce a weak homotopy equivalence $X\to X/\sim$ to the trajectory space \cite{BdH, CG}. Our result holds in fact with far more generality. In the case of a system of comparative retractions, the trajectory space is just the colimit of the diagram $X=X_0\to X_1\to \ldots $. The category $\aatop$ is closed under colimits. The analogous for limits is not true. In fact, Clader \cite{Bil, Cla, Thi} proved that any compact polyhedron is a strong deformation retract of a sequential limit of finite spaces. Our general result deals with filtered colimits where the maps in the diagram are weak equivalences. We prove the following

\begin{teomain}
	Let $F: \mathcal{C} \to \aatop$ be a functor from a filtered category to A-spaces such that every map in the diagram $F$ is a weak equivalence, then ${p_c:F(c)\to \text{colim }F}$ is a weak equivalence for every $c\in\mathcal{C}$.
\end{teomain} 

In order to prove Theorem \ref{main} we give a characterization of the minimal open sets of quotients of A-spaces, and prove a lifting property for maps to a filtered colimit.

Theorem \ref{main} is notably false in $\mathsf{Top}$, even if all the spaces $F(c)$ are polyhedra. Moreover, we will prove in Theorem \ref{graphs} that every locally finite $T_0$ space is a sequential colimit of graphs ($1$-dimensional polyhedra). In particular we obtain the following unexpected result.

\begin{intro2}
Every topological space has the weak homotopy type of a sequential colimit of graphs.
\end{intro2}

We use Theorem \ref{main} to prove that, although not every A-space has a core, it does have a \textit{weak core} (Theorem \ref{weakcore}). A weak core of $X$ is a space $Y$ such that there is a weak equivalence $X\to Y$ which is a quotient map, and where $Y$ is minimal with this property. Since the McCord map $|K|\to \x(K)$ is a quotient and a weak equivalence, we deduce the following 

\begin{intro1}
Every regular CW-complex admits a weak core.
\end{intro1}

In general, arbitrary topological spaces need not have weak cores (Example \ref{ejemplok}), and some spaces may have weak cores which are not Alexandroff (Example \ref{nota}).

We will see that in contrast to Stong's theory for finite spaces, weak cores of A-spaces are not unique in general (Example \ref{doswc}, Section \ref{sectionmany} and Section \ref{subdiv}).

Given any two spaces of the same weak homotopy type, there exists a common CW-approximation. From the results of Hardie and Vermeulen \cite{HV} it follows that if $X$ and $Y$ are finite spaces of the same weak homotopy type, there exist a finite space $Z$ and weak equivalences $X\leftarrow Z \to Y$. In \cite[Proposition 4.6.6]{Bar} it is proved that if $X$ and $Y$ are homotopy equivalent finite spaces, there exists a finite space $Z$ and quotient maps $X\leftarrow Z \to Y$ which are homotopy equivalences. In Section \ref{subdiv} we introduce a general notion of \textit{subdivision} of a $T_0$ A-space, extending barycentric subdivisions \cite{HV}, and prove an analogous of these results in our context: if $X$ and $Y$ are A-spaces of the same weak homotopy type, there exists an A-space $Z$ and quotient maps $X\leftarrow Z \to Y$ which are also weak equivalences.

\section{Preliminaries} \label{prelim}

Given an Alexandroff space $X$, its \textit{associated preorder} $\le$ on $X$ is given by $y\le x$ if $x$ lies in the closure $\overline{\{y\}}$ of $\{y\}$. The \textit{minimal open neighborhood} $U_x$ of a point $x$ is the intersection of all the open sets containing $x$. Then $y\le x$ if and only if $y\in U_x$. Conversely, if $\le$ is a preorder in a set $X$, a \textit{down-set} is a subset $U$ of $X$ such that $x\in U$ and $y\le x$ imply $y\in U$. The \textit{topology on $X$ associated to $\le$} is given by all the down-sets. These constructions furnish a correspondence between Alexandroff topologies on a set $X$ and preorders on the same set. A map between A-spaces is continuous if and only if it is order-preserving. Using this isomorphism, we will constantly think of A-spaces and preordered sets as the same thing. 

An A-space satisfies the $T_0$ separation axiom if and only if its preorder is an order. Every space $X$ deformation retracts to its Kolmogorov quotient $Q(X)$, the space obtained by identifying points with the same neighborhoods. Any section $s:Q(X)\to X$ is a homotopy inverse of the quotient. Two points $x,y$ in an A-space $X$ have the same neighborhoods if and only if $x\le y$ and $y\le x$ simultaneously. The Kolmogorov quotient $Q(X)$ of an A-space $X$ is an A-space and two classes $[x],[y]$ satisfy $[x]\le [y]$ if and only if $x\le y$. Any continuous map $f:X\to Y$ from a space to a $T_0$ space factorizes through $Q(X)$. In particular, any continuous map $f:X\to Y$ between A-spaces induces a continuous map $\overline{f}:Q(X)\to Q(Y)$ given by $\overline{f}([x])=[f(x)]$. Thus, $Q:\aatop\to \tatop$ is a functor from A-spaces to $T_0$ A-spaces.

If $X$ is a $T_0$ A-space, its \textit{associated simplicial complex} $\mathcal{K}(X)$ has as vertices the points of $X$ and as simplices the finite non-empty chains in $X$. The map $\mu:|\mathcal{K}(X)|\to X$ defined by $\mu(\alpha)=\text{min}(\text{sop}(\alpha))$ is called the \textit{McCord map}. This map is well-defined, continuous and, moreover, a weak homotopy equivalence, i.e. for every choice of base point $\alpha \in |\mathcal{K}(X)|$ and every $n\in\N_0$, ${\mu_*:\pi_n(|\mathcal{K}(X)|,\alpha)\to \pi_n(X,\mu(\alpha))}$ is a bijection. The proof by McCord that $\mu$ is a weak equivalence uses an important local to global criterion that we now recall in a simplified version.

\begin{teo}[McCord \cite{Mcc}] \label{mccord}
Let $X$ be an A-space and let $Y$ be any topological space. Let $f:Y\to X$ be a continuous map such that $f^{-1}(U_x)$ is homotopically trivial (that is, all its homotopy groups are trivial) for every $x\in X$. Then $f$ is a weak equivalence.
\end{teo}

If $f:X\to Y$ is a continuous map between $T_0$ A-spaces, there is a simplicial map $\kp (f):\kp (X)\to \kp (Y)$ which coincides with $f$ on vertices, and there is a commutative diagram

\[\begin{tikzcd}
		 {|\mathcal{K}(X)|} & {|\mathcal{K}(Y)|} \\
		{X} & {Y}
		\arrow["{|\mathcal{K}(f)|}", from=1-1, to=1-2]
		\arrow["{\mu}", from=1-1, to=2-1]
		\arrow["{\mu}", from=1-2, to=2-2]
		\arrow["f", from=2-1, to=2-2]
	\end{tikzcd}\]


If $K$ is a simplicial complex, the \textit{face poset} $\mathcal{X}(K)$ is the poset of simplices ordered by inclusion. In this case $\mathcal{K}\mathcal{X}(K)$ is the barycentric subdivision of $K$, so we have the McCord map $\mu:|K|=|\mathcal{K}\mathcal{X}(K)|\to \mathcal{X}(K)$. If $\varphi: K\to L$ is a simplicial map between simplicial complexes, there is a continuous map $\x (\varphi):\x(K)\to \x(L)$ defined by $\x(\varphi)(\sigma)=\varphi(\sigma)$. The square

\[\begin{tikzcd}
		 {|K|} & {|L|} \\
		{\x(K)} & {\x(L)}
		\arrow["{|\varphi|}", from=1-1, to=1-2]
		\arrow["\mu", from=1-1, to=2-1]
		\arrow["{\mu}", from=1-2, to=2-2]
		\arrow["\x(\varphi)", from=2-1, to=2-2]
	\end{tikzcd}\]
only commutes up to homotopy.

We say that a space $X$ is \textit{locally finite} if every point $x\in X$ has a finite neighborhood. Every locally finite space is Alexandroff.

\begin{obs} \label{cerrado}
	To see that $\aatop\subseteq \mathsf{Top}$ is closed under colimits we prove something slightly more general: If $Y$ is a topological space with the final topology induced by the family $\{f_i:X_i\to Y\}_{i\in \mathcal{I}}$ and every $X_i$ is Alexandroff, then $Y$ is Alexandroff. The claim then follows since a colimit of topological spaces is constructed as a quotient of a disjoint union, both of these stages defined through a final topology.
	
	Suppose $\{V_j\}_{j\in\mathcal{J}}$ is a family of open sets in $Y$, then $f_i^{-1}(V_j)\subseteq X_i$ is open for every $i$ and $j$. Since each $X_i$ is Alexandroff, ${f_i^{-1}(\underset{j\in\mathcal{J}}{\bigcap} V_j)=\underset{j\in\mathcal{J}}{\bigcap} f_i^{-1}(V_j)}$ is open for every $i$. Then $\underset{j\in\mathcal{J}}{\bigcap} V_j$ is open in $Y$.
\end{obs}

\begin{obs} \label{imagenfinita}
	A useful consequence of the work of McCord is that any element in a homotopy group of an Alexandroff space $X$ and any homotopy class of maps from a compact CW-complex to $X$, is represented by a map with finite image. 
	
	Let $X$ be a $T_0$ A-space, $\mu:|\mathcal{K}(X)|\to X$ its McCord map, and $n\in\N_0$. Since $\pi_n(|\mathcal{K}(X)|,x)\xrightarrow{\mu_*} \pi_n(X,x)$ is an isomorphism for every $x\in X$, there exists for every map $(S^n,s_0)\xrightarrow{\gamma}(X,x)$ another map $(S^n,s_0)\xrightarrow{\widetilde{\gamma}}(|\mathcal{K}(X)|,x)$ such that ${[\mu \widetilde{\gamma}]=[\gamma]\in\pi_n(X,x)}$.
	Since a compact subset of a simplicial complex intersects only finitely many open simplices, $\text{Im}(\mu \widetilde{\gamma})\subseteq X$ is finite by definition of $\mu$.
    
     If $L$ is a CW-complex, the McCord map induces a bijection $[L,|\mathcal{K}(X)|]\to [L,X]$ between unpointed homotopy classes. Thus, if $L$ is also compact, then every element in the codomain is represented by a map with finite image.

    If $X$ is now an A-space which is not necessarily $T_0$, the same statement holds as the Kolmogorov quotient is a strong deformation retract which contains any chosen point $x\in X$.
\end{obs}

\begin{ej} \label{ejemplononh}
    It is not true that every map from a compact topological space to an A-space is homotopic to a map with finite image. Let $X$ be the non-Hausdorff suspension of an infinite discrete space, that is the poset $X=D\cup\{a,b\}$, where $D$ is an infinite set and the order is given by $a>d, b>d$ for every $d\in D$. Then $X$ is compact and the identity $1_X$ is not homotopic to a map with finite image. This follows from an argument using fundamental groups (see Example \ref{ejl}).  
\end{ej}

If $f,g:X\to Y$ are continuous maps between A-spaces, we write $f\le g$ if $f(x)\le g(x)$ for every $x\in X$. In this case there is a homotopy from $f$ to $g$ which is identically $f$ for $t<1$, so, in fact, $f\simeq g \text{ rel } \{x: f(x)=g(x)\}$. A point $x$ in a $T_0$ A-space $X$ is said to be an up-beat point if the set $X_{>x}$ has a minimum $y$. In this case the retraction $r:X\to X\smallsetminus \{x\}$ which maps $x$ to $y$ is continuous, and moreover, a homotopy equivalence, as $1_X\le ir$. Here $i:X\smallsetminus \{x\}\to X$ denotes the inclusion. A down-beat point is defined dually, and its removal also gives a deformation retract. When $X$ is a finite space, $Q(X)$ is finite and $T_0$ and removing beat points (up or down) eventually leads to a space $Y$ without beat points. There is then a quotient $X\to Y$ and it can be proved that $Y$ has no proper quotients which are homotopy equivalent \cite{Sto}. The space $Y$ is called a core of $X$ and it is unique up to homeomorphism: it does not depend on the sequence of beat points chosen. 

Note that if $X$ is an A-space which has a point $x$ which is comparable to any other point, then $X$ is contractible. Indeed, if $f:X\to X$ maps the points $y\ge x$ to $x$ and fixes any other point, $f$ is continuous and $1_X\ge f\le c_x$, where $c_x$ is the constant map to $x$.

\section{Quotients and liftings}

Our next result describes the topology of any quotient of an A-space.

\begin{defi}
	Let $X$ be an A-space and $R$ an equivalence relation on $X$. Given $x,y\in X$, define an \textit{$R$-sequence} from $y$ to $x$ to be a finite sequence of elements $x_0,y_0,x_1,y_1,x_2,y_2\dots,$ $x_m,y_m \in X$ such that
	\begin{equation*}
		y=x_0Ry_0\le x_1Ry_1\le x_2Ry_2\le \dots \le x_mRy_m=x.
	\end{equation*}
	Define $R_x=\{y\in X: \text{an $R$-sequence from $y$ to $x$ exists}\}$. 
\end{defi}

If $X$ is a space and $R$ is an equivalence relation on $X$, we denote by $p:X\to \sfrac{X}{R}$ the canonical projection to the quotient, considered with the quotient topology. For $x\in X$, we denote by $[x]$ the class of $x$ in $\sfrac{X}{R}$.

\begin{lema} \label{cocientes}
    Let $X$ be an A-space and let $R$ be an equivalence relation on $X$. Then $R_x$ is open, saturated and $p^{-1}(U_{[x]})=R_x$ for every $x\in X$. Therefore, $U_{[x]}=p(R_x)$.
\end{lema}
\begin{proof}	
Take $y\in R_x$ and an $R$-sequence
${y=x_0Ry_0\le \dots \le x_mRy_m=x}$. If ${z\le y}$, then ${zRz\le yRy_0\le \dots \le x_mRy_m=x}$ is an $R$-sequence and ${z\in R_x}$. Since $U_y\subseteq R_x$ for every $y\in R_x$, $R_x$ is open.

To prove that $R_x$ is saturated it suffices to prove $p^{-1}(p(R_x))\subseteq R_x$. Let $y\in p^{-1}(p(R_x))$. Then there exists $z\in R_x$ such that $[z]=[y]$. As in the previous item, take an $R$-sequence from $z$ to $x$ and add the first term $yRz$ to see that $y\in R_x$.

Since $R_x$ is open and saturated, $p(R_x)$ is an open set that contains $[x]$. Therefore $U_{[x]}\subseteq p(R_x)$, which implies $p^{-1}(U_{[x]})\subseteq p^{-1}(p(R_x))=R_x$.

Conversely, if $y\in R_x$, there is an $R$-sequence ${y=x_0Ry_0\le x_1Ry_1\le \dots \le x_mRy_m=x}$. When passing to the quotient, $[y]=[x_0]\le [x_1]\le\dots\le [x]$ implies $[y]\in U_{[x]}$. Thus, $p(R_x)\subseteq U_{[x]}$ and then ${R_x=p^{-1}(p(R_x))\subseteq p^{-1}(U_{[x]})}$. 
\end{proof}

When given a functor $F:\mathcal{C}\to \aatop$ to A-spaces we will denote by $p_c$ the canonical map $F(c)\to \text{colim }F$ for every object $c\in\mathcal{C}$. Since $\aatop$ is closed by colimits in $\topp$, when $\mathcal{C}$ is filtered, the colimit $\text{colim } F$ can be constructed as the quotient of $\underset{c\in\mathcal{C}}{\coprod} F(c)$ by the equivalence relation $R$ given by: $x_1\in F(c_1)$ and $x_2\in F(c_2)$ are related if there exists $c\in\mathcal{C}$ and morphisms $c_1\xrightarrow{g_1} c$, $c_2\xrightarrow{g_2} c$ such that $F(g_1)(x_1)=F(g_2)(x_2)$. In this case, the map $p_c$ is the inclusion of $F(c)$ in the disjoint union, composed with the quotient map $p:\underset{c'\in\mathcal{C}}{\coprod} F(c') \to (\underset{c'\in\mathcal{C}}{\coprod} F(c'))/R$.

\begin{lema} \label{diez}
	Let $\mathcal{C}$ be a filtered category and $F:\mathcal{C}\to \aatop$ a functor. For every pair of elements $x,y\in \text{colim }F$ such that $x\le y$ there exists an object $c$ in $\mathcal{C}$ and elements $\widetilde{x},\widetilde{y} \in F(c)$ such that $\widetilde{x}\le \widetilde{y}$, $p_c(\widetilde{x})=x$ and $p_c(\widetilde{y})=y$. 
\end{lema}
\begin{proof}
	
	
	Let $a\in p^{-1}(x)$ and $b\in p^{-1}(y)$. Since $x\le y$, $x\in U_y$, and since ${p^{-1}(U_y)=R_b}$ by Lemma \ref{cocientes}, there exists an $R$-sequence from $a$ to $b$. Say, $a=x_0 R y_0\le x_1 R y_1\le x_2 R y_2\le \ldots \le x_m R y_m= b$. Suppose $m\ge 2$.

    Let $c_0,c_1,c_2 \in \mathcal{C}$ be objects such that $x_0\in F(c_0)$, $y_0,x_1\in F(c_1)$, $y_1,x_2\in F(c_2)$.  By definition of $R$ there exist $e_0,e_1\in \mathcal{C}$ and morphisms $c_0\xrightarrow{g_0}e_0$, $c_1\xrightarrow{h_0}e_0$, $c_1\xrightarrow{g_1}e_1$, $c_2\xrightarrow{h_1}e_1$ such that $F(g_i)(x_i)=F(h_i)(y_i)$ for $i=0,1$.

    Since $\mathcal{C}$ is filtered, there exists $c\in \mathcal{C}$ and morphisms $e_0\xrightarrow{g}c$, $e_1\xrightarrow{h}c$ such that $gh_0=hg_1$. Thus,

    $$a=x_0 R F(gg_0)(x_0)=F(gh_0)(y_0)\le F(gh_0)(x_1)=F(hg_1)(x_1)=F(hh_1)(y_1)\le$$ $$\le F(hh_1)(x_2) R x_2 R y_2.$$

    Thus, there is shorter $R$-sequence from $a$ to $b$.

    When $m=1$, the argument above shows that $\widetilde{x}=F(gg_0)(x_0)$ and $\widetilde{y}=F(hh_1)(y_1)$ satisfy the conditions of the statement.
\end{proof}



\begin{prop} \label{levantarx}
	Let $F:\mathcal{C}\to \aatop$ be a functor from a filtered category to A-spaces. If $X \subseteq \text{colim }F$ is a finite subspace, then there exists $c\in\mathcal{C}$ and $\widetilde{X}\subseteq F(c)$ such that $p_c(\widetilde{X})\subseteq X$ and $p_c \restr _{\widetilde{X}}: \widetilde{X} \to X$ is a homeomorphism.
	
\end{prop}
\begin{proof}
	
	
	
	Let $X\subseteq \text{colim }F$ be a finite subspace. By Lemma \ref{cocientes}, for each pair $x\le y$ in $X$ there exists $c_{x,y}\in \mathcal{C}$ and elements $a_{x,y}\le b_{x,y}$ in $F(c_{x,y})$ such that $p_{c_{x,y}}(a_{x,y})=x$ and $p_{c_{x,y}}(b_{x,y})=y$. 
    
    Since $X$ is finite and $\mathcal{C}$ is filtered, there exists $c\in\mathcal{C}$ and morphisms $c_{x,y}\xrightarrow{g_{x,y}} c$, such that $F(g_{x,y})(a_{x,y})=F(g_{z,x})(b_{z,x})$ for every $z\le x\le y$ in $X$.
		
	Define $\widetilde{X}=\{F(g_{x,x})(a_{x,x})\}_{x\in X} \subseteq F(c)$. It is easy to see that $p_c(\widetilde{X})\subseteq X$ and that $p_c\restr _{\widetilde{X}}:\widetilde{X}\to X$ is continuous and bijective. We prove that $(p_c\restr _{\widetilde{X}})^{-1}:X\to \widetilde{X}$ is also continuous. Let $x\le y$ in $X$. Then $$(p_c\restr _{\widetilde{X}})^{-1}(x)=F(g_{x,x})(a_{x,x})=F(g_{x,y})(a_{x,y})\le F(g_{x,y})(b_{x,y})=F(g_{y,y})(a_{y,y})=(p_c\restr _{\widetilde{X}})^{-1}(y).$$

\end{proof}

\begin{coro} \label{tecero}
	If $F:\mathcal{C}\to \aatop$ is a functor from a filtered category and $F(c)$ is $T_0$ for every $c\in \mathcal{C}$, then $\text{colim }F$ is $T_0$. 
\end{coro}
\begin{proof} 	
	By Proposition \ref{levantarx}, any two-element subspace of $\text{colim }F$ is homeomorphic to a subspace of some $F(c)$, which is $T_0$. Thus, $\text{colim }F$ is $T_0$. 
\end{proof}

The previous result is not true in general when $\mathcal{C}$ is not filtered. The $T_0$ space of three points $a<b<c$ has a non-$T_0$ quotient obtained by identifying $a$ and $c$. Thus, there is a pushout of $T_0$ A-spaces which is not $T_0$.

\begin{ej} \label{ejt0}
A filtered colimit of $T_0$ spaces in general is not $T_0$. Consider $S^1\subseteq \C$ and the map $g:S^1\to S^1$ given by 
	\begin{equation*}
		g(e^{\theta i})=\begin{cases}
			e^{2 \theta i} \quad \text{if }0\le \theta\le \pi/2\\
            -1 \quad \text{if } \pi/2\le \theta\le \pi\\
            e^{(2\theta -\pi) i} \quad \text{if }\pi\le \theta\le 3\pi/2\\
			1 \quad \text{if }3\pi/2\le \theta\le 2\pi.
		\end{cases}
	\end{equation*}
	Essentially, $g$ stretches the first quarter of $S^1$ so it fills the first half. It collapses the second quarter to $-1$, it stretches the third quarter to the second half, and collapses the fourth quarter to $1$.
	
	Define the functor $F:\N \to \mathsf{Top}$ by $F(n)=S^1$ for every $n\in \N$ and $F(n\to n+1)=g$. All the points $e^{\theta i}$ with $0<\theta \le \pi$ eventually reach $-1$ and all the other eventually reach $1$. Thus $\text{colim }F$ has just two points, and it is easy to see that none of them is open. Thus, $\text{colim }F$ is not $T_0$.
\end{ej} 




From Proposition \ref{levantarx} we immediately deduce:

\begin{coro} \label{lift}
	Let $F:\mathcal{C}\to \aatop$ be a functor from a filtered category to A-spaces. If $f:Y\to \text{colim }F$ is a continuous map from any topological space $Y$ with finite image, then there exists $c\in\mathcal{C}$ and $\widetilde{f}:Y\to F(c)$ such that $p_c \widetilde{f}=f$.
	\[\begin{tikzcd}
		& {F(c)} \\
		Y & {\text{colim }F}
		\arrow["{p_c}", from=1-2, to=2-2]
		\arrow["{\widetilde{f}}", from=2-1, to=1-2]
		\arrow["f"', from=2-1, to=2-2]
	\end{tikzcd}\]
	Moreover, the map $\widetilde{f}$ can be taken in such a way that $\text{Im}\widetilde{f}\xrightarrow{p_c}\text{Im}f$ is a homeomorphism.
\end{coro}

Corollary \ref{lift} can be used to show that finite spaces are compact objects in the category $\aatop$. The converse is also true as A-spaces are filtered colimits of their finite subspaces. 

\section{Weak equivalences and colimits}

If $F:\mathcal{C}\to \mathsf{Top}$ is a functor from a filtered category such that every map in the diagram $F$ is a weak equivalence, then it is not necessarily true that the maps $p_c:F(c)\to\text{colim }F$ are weak equivalences, even if the topological spaces $F(c)$ are polyhedra. This is already apparent from Example \ref{ejt0}, since the map $g$ defined there is a homotopy equivalence, but the colimit is the indiscrete space of two points, which is contractible. A similar example can be given by replacing $g$ by a map which stretches the first half of $S^1$ to all of $S^1$ while collapses the second half to $1\in S^1$. This is also a homotopy equivalence and the colimit is the singleton. Our main result in this section is that a filtered colimit in Alexandroff spaces preserves weak equivalences.

\begin{teo} \label{main}
	Let $F: \mathcal{C} \to \aatop$ be a functor from a filtered category to A-spaces such that every map in the diagram $F$ is a weak equivalence, then ${p_c:F(c)\to \text{colim }F}$ is a weak equivalence for every $c\in\mathcal{C}$.
\end{teo}
\begin{proof}	
We first observe that all spaces $F(c)$ may be assumed to be $T_0$. Indeed, it is easy to check that the Kolmogorov quotient functor $Q: \aatop \to \tatop$ is a left adjoint to the inclusion $i:\tatop \to\aatop$, so in particular it preserves colimits. By Corollary \ref{tecero} the filtered colimit $\text{colim }QF \in \tatop$ coincides with the colimit $\text{colim } iQF \in \aatop$. By the naturality of the projections to Kolmogorov quotients we have the following commutative diagram	of topological spaces \[\begin{tikzcd}
		{F(c)} &&& {Q(F(c))} \\
		{\text{colim }F} \\
		{Q(\text{colim }F)} &&& {\text{colim }QF}
		\arrow[from=1-1, to=1-4]
		\arrow["{p_c}", from=1-1, to=2-1]
		\arrow["{Q(p_c)}"', from=1-4, to=3-1]
		\arrow["p_{QF(c)}", from=1-4, to=3-4]
		\arrow[from=2-1, to=3-1]
		\arrow["=", from=3-4, to=3-1]
	\end{tikzcd}\]
	The projections are homotopy equivalences, so if $p_{QF(c)}:QF(c)=iQF(c)\to \text{colim }iQF=\text{colim }QF$ is a weak equivalence, then $p_c$ is a weak equivalence. We assume then that each space $F(c)$ is $T_0$.
    
	Consider $c\in \mathcal{C}$. We will prove that $(p_c)_*:\pi_k(F(c),x)\to \pi_k(\text{colim }F,p_c(x))$ is a monomorphism and an epimorphism for every $k\in \N_0$ and $x\in F(c)$.

	We prove first that it is an epimorphism. Take $[\gamma]\in\pi_k(\text{colim }F,p_c(x))$ with $\gamma: (S^k,s_0)\to (\text{colim }F,p_c(x))$. We can assume $\gamma$ has finite image by Remark \ref{imagenfinita}. Then, by Corollary \ref{lift}, there exist $c'\in\mathcal{C}$ and a map ${\widetilde{\gamma}:S^k\to F(c')}$ such that $p_{c'} \widetilde{\gamma}=\gamma$.

	Then $x'=\widetilde{\gamma}(s_0)\in F(c')$ is in the same equivalence class as $x\in F(c)$. Since $\mathcal{C}$ is filtered, there exist $d\in\mathcal{C}$ and morphisms $c\xrightarrow{g_c}d,c'\xrightarrow{g_{c'}}d$ such that ${F(g_c)(x)=F(g_{c'})(x')=:y}$. We have a commutative diagram
	
	\[\begin{tikzcd}
		{\pi_k(F(c),x)} & {\pi_k(F(d),y)} & {\pi_k(F(c'),x')} \\
		\\
		& {\pi_k(\text{colim }F, p_c(x))}
		\arrow["{F(g_c)_*}", from=1-1, to=1-2]
		\arrow["{(p_c)_*}"{description}, from=1-1, to=3-2]
		\arrow["{(p_d)_*}"', from=1-2, to=3-2]
		\arrow["{F(g_{c'})_*}"', from=1-3, to=1-2]
		\arrow["{(p_{c'})_*}"{description}, from=1-3, to=3-2]
	\end{tikzcd}\]
	
	Since $F(g_c)_*$ is an isomorphism by hypothesis,     $[\gamma]=(p_{c'})_*([\widetilde{\gamma}])=$ $(p_d)_* F(g_{c'})_*([\widetilde{\gamma}])\in$ $\text{Im}((p_{d})_*)=\text{Im}((p_d)_* F(g_c)_*)=\text{Im}((p_c)_*)$. We conclude that $(p_c)_*$ is an epimorphism.

	Now we show that $(p_c)_*$ is a monomorphism. Since $\mu: |\mathcal{K}(F(c))|\to F(c)$ is a weak equivalence, any class in $\pi_k(F(c),x)$ can be represented by $\mu \gamma$ for some $\gamma: (S^k,s_0)\to (|\mathcal{K}(F(c))|,x)$.
	
	Suppose $(p_c)_*[\mu \gamma]=0\in\pi_k(\text{colim }F,p_c(x))$. Consider the following diagram.

    \[\begin{tikzcd}
		& {|\mathcal{K}(F(c))|} && {|\mathcal{K}(\text{colim }F)|} \\
		{S^k} & {F(c)} && {\text{colim }F}
		\arrow["{|\mathcal{K}(p_c)|}", from=1-2, to=1-4]
		\arrow["\mu", from=1-2, to=2-2]
		\arrow["{\mu'}", from=1-4, to=2-4]
		\arrow["\gamma", from=2-1, to=1-2]
		\arrow["{\mu \gamma}"', from=2-1, to=2-2]
		\arrow["{p_c}", from=2-2, to=2-4]
	\end{tikzcd}\]
where $\mu, \mu'$ denote the McCord maps.
	
	Then ${[\mu' |\mathcal{K}(p_c)| \gamma]=[p_c \mu \gamma]=0}$, and, since $\mu'$ is a weak equivalence, $[|\mathcal{K}(p_c)| \gamma] =0\in\pi_k(|\mathcal{K}(\text{colim }F)|, p_c(x))$. In particular there is a null-homotopy ${G: S^k\times I\to |\mathcal{K}(\text{colim }F)}|$, $G:|\mathcal{K}(p_c)| \gamma\simeq \text{const}_{p_c(x)}$. Define ${H=\mu' G:S^k\times I\to  \text{colim }F}$. Then $H$ is a homotopy between $p_c \mu \gamma$ and the constant map $\text{const}_{p_c(x)}$ to $p_c(x)$. Moreover, since $S^k\times I$ is compact, the image of $G$ intersects only finitely many open simplices, so $H$ has finite image. By Corollary \ref{lift} there exist $c'\in\mathcal{C}$ and a lift $\widetilde{H}:S^k\times I\to F(c')$ in such a way that $p_{c'}$ induces a homeomorphism between the images of $\widetilde{H}$ and $H$. 
	
	Then $\widetilde{H}$ is a homotopy between $\widetilde{H}_0$ and a constant map, and $p_{c'} \widetilde{H}_0=p_c \mu \gamma$. Since both $\widetilde{H}_0$ and $\mu \gamma$ have finite image and the same projection to the filtered colimit, there exist $d\in\mathcal{C}$ and morphisms $c\xrightarrow{g_c}d$, $c'\xrightarrow{g_{c'}}d$ such that ${F(g_c) \mu \gamma=F(g_{c'}) \widetilde{H}_0}$.
	
	Finally, since $\widetilde{H}_0$ is null-homotopic, $[\widetilde{H}_0]=0\in\pi_k(F(c'),\widetilde{H}_0(s_0))$, and then $0=[F(g_{c'}) \widetilde{H}_0]=[F(g_c) \mu \gamma]\in\pi_k(F(d),F(g_{c})(x))$. Since $F(g_{c})$ is a weak equivalence, this implies that ${[\mu \gamma]=0\in\pi_k(F(c),x)}$, as we wanted to prove.
\end{proof}

\begin{obs}
    Note that the proof of Theorem \ref{main} shows in fact that if $F: \mathcal{C} \to \aatop_*$ is a functor from a filtered category $\mathcal{C}$ to pointed A-spaces and every map in the diagram $F$ induces an isomorphism in $\pi_k$ for a fixed $k\ge 0$, then $(p_c)_*:\pi_k(F(c))\to \pi_k(\text{colim }F)$ is an isomorphism for every $c\in\mathcal{C}$.
\end{obs}

We comment on possible alternative proofs of Theorem \ref{main}. A different approach could have been given using that the analogous result holds in the category of simplicial complexes with simplicial maps (which can be proved by taking simplicial approximations of maps from spheres), and the fact that both functors $\mathcal{K}$ and $\mathcal{X}$ preserve filtered colimits. Conversely, the result for simplicial complexes can be directly deduced from Theorem \ref{main}.

A second alternative to prove that $(p_c)_*:\pi_k(F(c))\to \pi_k(\text{colim }F)$ is a monomorphism, avoiding partially the associated simplicial complexes, would be to try to lift a homotopy $H: S^k \times I \to \text{colim } F$ to some stage $F(c')$. However, $H$ may have infinite image, even when $H_0=p_c\mu \gamma$ and $H_1=\text{const}_{p_c(x)} $ have finite image. Here we could invoke Kong and Khalimsky results \cite{KK} or Wofsey's \cite{Wof} to lift $H_0$ and $H_1$ to $|\mathcal{K}(\text{colim } F)|$ and then project a homotopy between the lifts to obtain a homotopy with finite image between $H_0$ and $H_1$.

Yet another approach would be to reduce the proof of Theorem \ref{main} to locally finite spaces. This can be achieved by using that both $\mathcal{K}$ and $\mathcal{X}$ preserve filtered colimits. Then, the image of a homotopy $S^k\times I\to \text{colim F}$ is compact in a locally finite space, and thus, finite.

These ideas for obtaining homotopies with finite image can be used to prove that for any functor $F:\mathcal{C}\to \aatop_*$ from a filtered category to pointed A-spaces we have that $\text{colim}(\pi_k F)$ is isomorphic to $\pi_k(\text{colim }F)$.

\section{Alexandroff spaces as colimits of graphs}

By generalizing the ideas of Example \ref{ejt0}, and to understand to what extent Theorem \ref{main} fails for general spaces, we will prove that any topological space is weakly equivalent to a sequential colimit of graphs (CW-complexes of dimension 1), with all the maps in the diagram being homotopy equivalences.

Let $X$ be a poset. Denote by $\prec$ the covering relation on $X$, i.e. $x\prec y$ if $x< y$ and there is no $x<z<y$. In a locally finite $T_0$ space the transitive and reflexive closure of $\prec$ is the order associated to the topology.

\begin{defi}
    Given a locally finite $T_0$ space $X$, define $G(X)$ as the underlying graph of the Hasse diagram of $X$. In other words, it is the CW-complex with one 0-dimensional cell for each point in $X$ and one 1-dimensional cell attached to $x,y$ for each pair $x\prec y$. Equivalently $G(X)$ is the disjoint union $$(\coprod\limits_{x\prec y} I\times \{(x,y)\}) {\scriptstyle\coprod} (\coprod\limits_{x\in X} \{x\}) $$ with identifications $(0,(x,y))\sim x$ and $(1,(x,y))\sim y$. 
\end{defi}

\begin{ej}
For each $n\ge 0$ there is a space $\mathbb{S}^n$ of $2n+2$ points weak equivalent to $S^n$: $\mathbb{S}^n=\{0,1,\ldots, n\}\times \{1,-1\}$, with $(i,j)<(i',j')$ if $i<i'$. If $n\ge 1$, the graph $G(\mathbb{S}^n)$ has 4 vertices of degree 2 and $2n-2$ of degree 4. For $n\ge 2$, $G(\mathbb{S}^n)$ is homeomorphic to a regular graph of degree 4 and $2n-2$ vertices.
\begin{figure}[ht] 
\begin{center}
\includegraphics[scale=0.25]{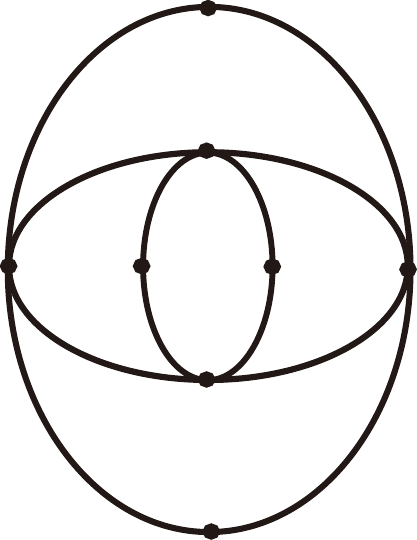}
\caption{The graph $G(\mathbb{S}^3)$ of the finite model of $S^3$ with 8 points.}\label{ojo}
\end{center}
\end{figure}
\end{ej}

\begin{teo} \label{graphs}
    Every locally finite $T_0$ space is a sequential colimit of graphs.
\end{teo}
\begin{proof}
    Let $X$ be a locally finite $T_0$ space and let $f:I\to I$ be defined by $$f(s)=\begin{cases}
			0  \quad \text{if }s\le \frac{1}{2} \\
			2s-1 \quad \text{if }s\ge \frac{1}{2}.
	\end{cases}$$
     This map induces a continuous map $\overline{f}:G(X)\to G(X)$ which fixes $[x]$ for every $x\in X$ and where $\overline{f}([(s,(x,y))])=[(f(s),(x,y))]$. Since $f\simeq 1_I \text{ rel } \{0,1\}$, $\overline{f}\simeq 1_{G(X)}$. In particular, $\overline{f}$ is a homotopy equivalence. 

     Consider the colimit $Y$ of the following diagram:
\[\begin{tikzcd}
	{G(X)} & {G(X)} & {G(X)} & \dots
	\arrow["\overline{f}", from=1-1, to=1-2]
	\arrow["\overline{f}", from=1-2, to=1-3]
	\arrow["\overline{f}", from=1-3, to=1-4]
\end{tikzcd}\]

The trajectory of a point $[(s,(x,y))]$ with $s<1$, in any copy of $G(X)$, eventually reaches a fixed point $[x]$. Thus $h:X\to Y$ defined by $h(x)=[([x],0)]$ is a bijection, where $[([x],0)]$ denotes the class of $[x]$ in the first copy of $G(X)$. Let $U\subseteq Y$ be open, and suppose $x\in h^{-1}(U)$, that is $[([x],0)]\in U$. Then $p_0^{-1}(U)$ is an open subset of $G(X)$ which contains $[x]$. If $y\prec x$, then $p_0^{-1}(U)$ contains a point $[(s,(x,y))]$ for some $s<1$. Thus, $[([y],0)]=p_0([(s,(x,y))])\in U$, and then $y\in h^{-1}(U)$. Since $X$ is locally finite, this implies that $U_x\subseteq h^{-1}(U)$, so $h^{-1}(U)$ is open. This proves that $h$ is continuous.

We prove now that $h$ is open. Let $U\subseteq X$ be open. In order to see that $h(U)\subseteq Y$ is open, we prove that $p_n^{-1}(h(U))\subseteq G(X)$ is open for every $n\ge 0$. If a point $[(s,(x,y))]$ with $0<s<1$ lies in $p_n^{-1}(h(U))$, then all the points $[(t,(x,y))]$ with $0<t<1$ project to the same point of $Y$ and, thus, they are also in $p_n^{-1}(h(U))$. Therefore, $p_n^{-1}(h(U))$ is a neighborhood of $[(s,(x,y))]$. If a point $[x]$ is in $p_n^{-1}(h(U))$, then all the points $[(s,(x,y))]$ with $s<1$ project to the same point, and then they are also in $p_n^{-1}(h(U))$. The points $[(s,(y,x))]$ with $s<1$ project to $[([y],0)]$, which is contained in $h(U)$ as $U\subseteq X$ is open and $y\le x$. Thus $[(s,(y,x))]\in p_n^{-1}(h(U))$ for any $s\in I$ and $y\prec x$. This proves that $p_n^{-1}(h(U))$ is a neighborhood of $[x]$. This concludes the proof that $h$ is open, and therefore, a homeomorphism.
\end{proof}

\begin{obs}
    By adding loops to the graph $G(Q(X))$ and using the ideas in Example \ref{ejt0} one can extend Theorem \ref{graphs} and prove that every locally finite space $X$ (not necessarily $T_0$) is a sequential colimit of graphs.
\end{obs}

Since every space is weak homotopy equivalent to a polyhedron and the functor $\mathcal{X}$ assigns to every simplicial complex a weakly equivalent locally finite and $T_0$ space, we have the desired result.

\begin{coro} \label{maincoro2}
	Every topological space has the weak homotopy type of a sequential colimit of graphs.
\end{coro}

\section{Weak Cores}

In this section we introduce the notion of a weak core and use Theorem \ref{main} to prove its existence for A-spaces and regular CW-complexes. 

\begin{defi}
	Let $X$ be a topological space. A \textit{weak core} of $X$ is a topological space $Y$ together with a quotient map $X\to Y$ which is also a weak equivalence and such that any quotient map $Y\to Z$ which is a weak equivalence is a homeomorphism. We will also say that the space $Y$ alone is a weak core of $X$, without explicit mention to the map. 
\end{defi}

Quotient maps are convenient for different reasons. On the one hand, if $\mathcal{C}$ is a category with initial object $0$, and $F:\mathcal{C}\to \mathsf{Top}$ is a functor such that every map in the diagram is a quotient, then $\text{colim } F$ is a quotient of $X=F(0)$. Namely, it is homeomorphic to $X/\sim$ where $x\sim x'$ for $x,x' \in X$ if there exists $c\in \mathcal{C}$ such that $F(0\to c)(x)=F(0\to c)(x')$.
On the other hand the category of A-spaces is closed by quotients (Remark \ref{cerrado}).

Finally, the notion of core, that we generalize here, can be stated in terms of quotient maps, as was done in the Section \ref{intro}. 

\begin{teo} \label{weakcore}
	Let $X$ be an A-space. Then there exists a weak core of $X$.
\end{teo}
\begin{proof}
	Consider the set $E$ of equivalence relations $r\subseteq X\times X$ over $X$ such that $X\to X/r$ is a weak equivalence. The inclusion makes $E$ an ordered set.

    Let $\mathcal{C}$ be a chain in $E$. Let $\overline{\mathcal{C}}=\mathcal{C}\cup \{1_X\}$. We see $\overline{\mathcal{C}}$ as a category in the usual way, namely there is a unique morphism $r\to r'$ if $r\subseteq r'$ and there is no morphism otherwise. We define a functor $F:\overline{\mathcal{C}}\to \aatop$ which maps an equivalence relation $r$ to the space $X/r$, which is an A-space by Remark \ref{cerrado}. If $r,r'\in \overline{\mathcal{C}}$ and $r\subseteq r'$, define $F(r\to r')$ to be the map $X/r\to X/r'$ that maps the class of $x$ to the class of $x$. Note that this map is a weak equivalence as both $X\to X/r$ and $X\to X/r'$ are weak equivalences. Since $\overline{\mathcal{C}}$ is a filtered category, $X=F(1_X)\to \text{colim }F$ is a weak equivalence by Theorem \ref{main}. In this case $\text{colim }F$ is homeomorphic to $X/R$ where $R=\bigcup\limits_{r\in \overline{\mathcal{C}}} r$. Thus, $R\in E$ is an upper bound of $\mathcal{C}$.

    By Zorn's lemma, there exists a maximal element $r$ in $E$. We claim that $q:X\to X/r$ is a weak core of $X$. If $q':X/r\to Y$ is a quotient map which is a weak equivalence, define a relation $r'$ over $X$ by $x_1r'x_2$ if $q'q(x_1)=q'q(x_2)$. Then $Y$ is homeomorphic to $X/r'$, $r'\in E$ and $r\subseteq r'$. By maximality of $r$, $r'=r$, so $q'$ is a homeomorphism. 
\end{proof}

In order to generalize this result to polyhedra we only need to observe that the McCord map $|K|\xrightarrow{\mu} \mathcal{X}(K)$ is a quotient map besides a weak equivalence for every simplicial complex $K$.

\begin{prop}
    Let $X$ be a $T_0$ A-space. Then $\mu :|\mathcal{K}(X)|\to X$ is a quotient map.
\end{prop}
\begin{proof}
    We already know that $\mu$ is continuous and surjective.
    Let $U\subseteq X$ be a subset such that $\mu^{-1}(U)\subseteq |\mathcal{K}(X)|$ is open. We want to prove that $U$ is open. Let $x\in U$ and let $y<x$. Then $x$ is a vertex of $\mathcal{K}(X)$ and $\mu(x)=\min \text{sop}(x)=\min \{x\}=x\in U$, so $x\in \mu^{-1}(U)$. Since $y<x$, $xy$ is a 1-dimensional simplex of $\mathcal{K}(X)$. Since $\mu^{-1}(U)$ is open, there is a point $\alpha$ other than $x$ in $\mu^{-1}(U)\cap |xy|$. Thus $y=\min \text{sop}(\alpha)=\mu(\alpha)\in U$. This proves that $U$ is open.  
\end{proof}


If $X$ is a regular CW-complex, it admits a triangulation $K$, and since $\mu:K\to \mathcal{X}(K)$ is a quotient, any weak core of $\mathcal{X}(K)$ is a weak core of $X$. Therefore we obtain the following

\begin{coro} \label{maincoro}
	Every regular CW-complex admits a weak core.
\end{coro}

In the following sections we analyze examples and properties of weak cores. We finish this section by observing that a space might not have weak cores, or have more than one, and that weak cores are not necessarily A-spaces. In Section \ref{sectionmany} we will prove that an A-space may have infinitely many weak cores and in Section \ref{subdiv} we will prove that for finitely many weak equivalent A-spaces $X_1,X_2,\ldots, X_n$, each of them a weak core (of itself), there is always an A-space $X$ for which each $X_i$ is a weak core.

\begin{ej} \label{doswc}
Let $X$ be the finite $T_0$ space in Figure \ref{figseis}.

\begin{figure}[ht] 
\begin{center}
\includegraphics[scale=0.4]{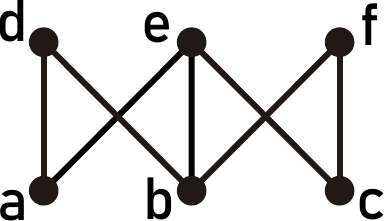}
\caption{A finite space with two weak cores.}\label{figseis}
\end{center}
\end{figure}

The quotient which identifies $d$ and $f$ gives a weak equivalence (by McCord's Theorem \ref{mccord}) to a space $Y$. Identifying $a$ and $c$ instead, gives a second weak equivalence to space $Z$. The fact that both $Y$ and $Z$ are weak cores of $X$ follows from the fact that any two points of any of them are part of a loop with image of cardinality 4 which is non-trivial in $\pi_1$. On the other hand, any loop with image of cardinality at most $3$ must be trivial in $\pi_1$. Thus, any proper quotient of $Y$ or $Z$ would induce a non-injective morphism in fundamental groups. This idea will be extended in Section \ref{sectionmany}. 

The spaces $Y$ and $Z$ are dual. Thus, the barycentric subdivision $Y'=\x \kp (Y)$ coincides with $Z'$. By \cite{HV} there is a weak equivalence $\zeta:Y'\to Y$, which in fact is also a quotient map. Thus $Y$ and $Z$ are also weak cores of $Y'$. Much more about subdivisions will be said in Section \ref{subdiv}. 
\end{ej}



\begin{ej} \label{ejemplok}
    There are spaces which do not have weak cores. Let $K=\{0\}\cup \{\frac{1}{n}\}_{n\in \N}\subseteq \R$ with the usual topology. Let $X=K\cup \{m\}$. The open subsets of $X$ are those of $K$, $X\smallsetminus \{0\}$ and $X$. 
    We claim that $X$ is not path-connected. Indeed, assume $\gamma :I\to X$ is a path from $0$ to $m$. Since $\{0\}\subseteq X$ is closed, there is a maximum $t_0$ of $\gamma^{-1}(0)$, and since $\{m\}$ is closed, there is a minimum $t_m$ of $\gamma^{-1}(m)\cap [t_0,1]$. Then $\gamma((t_0,t_m))$ is a connected subspace of $\{\frac{1}{n}\}_{n\in \N}$, so it is just one point. Thus, $\gamma$ is not continuous at $t_0$, a contradiction. 
    Note that the subspace $X\smallsetminus \{0\}$ of $X$ is path-connected, and in fact, contractible, since it is an A-space with maximum $m$.

    Suppose $q:X\to Z$ is a weak core of $X$. Note that $q(0)\neq q(x)$ for any $x\neq 0$, since otherwise $q_*:\pi_0(X)\to \pi_0(Z)$ would not be injective. 

    Since $X\smallsetminus \{0\}$ is a saturated open set, $q\restr_{X\smallsetminus \{0\}}: X\smallsetminus \{0\}\to Z\smallsetminus \{q(0)\}$ is also a quotient, and in particular $Z\smallsetminus \{q(0)\}$ is an A-space.
    
    Suppose there exists $n\in \N$ such that $q(\frac{1}{n})\neq q(m)$. Since $\frac{1}{n}<m$ in $X\smallsetminus \{0\}$, then $q(\frac{1}{n})\le q(m)$ in $Z\smallsetminus \{q(0)\}$. Assume there exists $k\in \N$ such that $q(\frac{1}{k})> q(\frac{1}{n})$ (meaning that $q(\frac{1}{k})\ge q(\frac{1}{n})$ and $q(\frac{1}{k})\neq q(\frac{1}{n})$). If $q(\frac{1}{k})\neq q(m)$, then $q^{-1}(q(\frac{1}{k}))\subseteq X\smallsetminus \{0\}$ does not contain the point $m$, so it is open, and then $\{q(\frac{1}{k})\}$ is open, which contradicts the fact that $q(\frac{1}{n})\in U_{q(\frac{1}{k})}$. This proves that the only point in $Z\smallsetminus \{q(0)\}$ greater than (and different from) $q(\frac{1}{n})$ is $q(m)$. Thus $Z\smallsetminus \{q(0)\}$ strong deformation retracts to $Z\smallsetminus \{q(0),q(\frac{1}{n})\}$. We show that this induces a strong deformation retraction from $Z$ to $Z\smallsetminus \{q(\frac{1}{n})\}$. To prove this note first that $q(\frac{1}{k})=q(\frac{1}{n})$ only for finitely many values of $k$. Otherwise, if $U$ is an open neighborhood of $q(0)$, $q^{-1}(U)$ is an open neighborhood of $0$, so it contains a point $\frac{1}{k}$ with $q(\frac{1}{k})=q(\frac{1}{n})$. This implies that $q(\frac{1}{n})\in U$. Since this is true for every open neighborhood $U$ of $q(0)$, $q(\frac{1}{n})$ and $q(0)$ are in the same path component. In fact, there is a path which is $q(\frac{1}{n})$ for $t<1$ and $q(0)$ for $t=1$. This contradicts the fact that $q$ induces an injection on $\pi_0$. Therefore, there exists $l\in \N$ such that $q(\frac{1}{k})\neq q(\frac{1}{n})$ for every $k\ge l$. Since $X\smallsetminus \{0\}$ and $\{0\}\cup \{\frac{1}{l}, \frac{1}{l+1}, \ldots\}$ are open and cover $X$, $X$ has the final topology with respect to the inclusions of these subspaces, and since $q$ is a quotient, $Z$ has the final topology with respect to the subspaces $Z\smallsetminus \{q(0)\}$ and $Y=q(\{0\}\cup \{\frac{1}{l}, \frac{1}{l+1}, \ldots\})$. Now, the retraction $r:Z\smallsetminus \{q(0)\}\to Z\smallsetminus \{q(0),q(\frac{1}{n})\}$ and the homotopy $H:(Z\smallsetminus \{q(0)\})\times I\to Z\smallsetminus \{q(0)\}$ between $ir$ and $1_{Z\smallsetminus \{q(0)\}}$ relative to $Z\smallsetminus \{q(0),q(\frac{1}{n})\}$ extend to a retraction $\widetilde{r}:Z\to Z\smallsetminus \{q(\frac{1}{n})\}$, which is continuous because it is on $Z\smallsetminus \{q(0)\}$ and $Y$, and to a homotopy $\widetilde{H}:Z\times I\to Z$ between $\widetilde{i}\widetilde{r}$ and $1_Z$, which is continuous because it is on $(Z\smallsetminus \{q(0)\})\times I$ and $Y\times I$. Here $i:Z\smallsetminus \{q(0),q(\frac{1}{n})\}\to Z\smallsetminus \{q(0)\}$ and $\widetilde{i}: Z\smallsetminus \{q(\frac{1}{n})\} \to Z$ denote the inclusions.

    Finally, since $\widetilde{r}$ is a quotient map, a homotopy equivalence and it is not a homeomorphism, $Z$ is not a weak core of $X$, a contradiction. We conclude that $q(\frac{1}{n})=q(m)$ for every $n\in \N$. Then $Z$ has only two points: $q(m)$, which is open, and $q(0)$, which is not. Thus $Z$ is the Sierpinski space, which is path-connected. This contradicts again the fact that $q_*:\pi_0(X)\to \pi_0(Z)$ is injective. 
\end{ej}

\begin{ej} \label{nota}
    There exist spaces which have weak cores that are not A-spaces. We start with the space $K=\{0\}\cup \{\frac{1}{n}\}_{n\in \N}$ of Example \ref{ejemplok}, but now we add two points instead of just one. Define $X=K\cup \{m,m'\}$. The open sets of $X$ are those of $K$ together with $X\smallsetminus \{0,m'\}$, $X\smallsetminus \{m'\}$, $X\smallsetminus \{0,m\}$, $X\smallsetminus \{m\}$, $X\smallsetminus \{0\}$ and $X$. Thus, $X\smallsetminus \{0\}$ is an A-space with two maximal points and countably many minimal points. It is in fact the non-Hausdorff suspension in Example \ref{ejemplononh}. An argument as in Example \ref{ejemplok} shows that $X$ has two path components, $\{0\}$ and $X\smallsetminus \{0\}$. 

    Suppose $q:X\to Z$ is a quotient map which is a weak equivalence. Then the restriction $q\restr_{X\smallsetminus \{0\}}: X\smallsetminus \{0\}\to Z\smallsetminus \{q(0)\}$ is a quotient and a weak equivalence. As explained in Example \ref{doswc}, $X\smallsetminus \{0\}$ is a weak core (of itself), as any two points are part of a loop with image of cardinality 4 that is non-trivial in $\pi_1$. Thus $q\restr_{X\smallsetminus \{0\}}$ is injective, and then so is $q$. Therefore $q$ is a homeomorphism. This proves that $X$ is a weak core of itself. Note that $X$ is not an A-space as the intersection $\{0\}$ of the open sets $V_n=\{0\}\cup \{\frac{1}{k} : k\ge n\}$ is not open. 

    The space $X$ is also a weak core of a subspace of the plane, see Figure \ref{figrombo}.

    \begin{figure}[ht] 
\begin{center}
\includegraphics[scale=0.4]{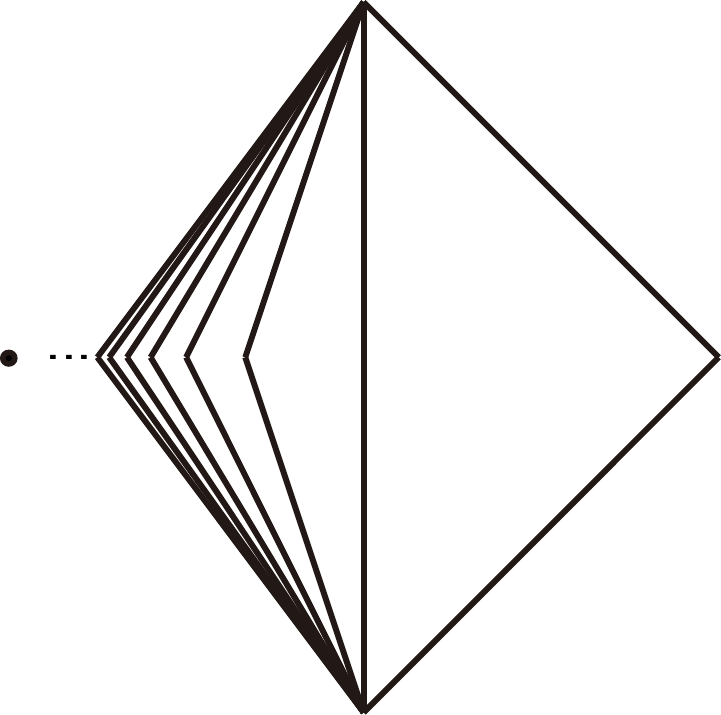}
\caption{A subspace of $\R^2$ with a weak core which is not an A-space.}\label{figrombo}
\end{center}
\end{figure}
\end{ej}

\section{Dismantlings of posets}

In the context of order theory, a dismantling of a poset $P$ is the procedure of applying successive retractions of a suitable class to find a minimal retract which admits no non-trivial such retraction. The class of retractions chosen depends on the property that is to be studied. The minimal retract possesses such a property if and only if the original poset $P$ has it.  The aim is to reduce the problem to a simpler combinatorial structure.

For a finite poset, Stong \cite{Sto} studied its homotopy type by applying retractions that remove beat points (see Section \ref{prelim}). The same idea was applied in order theory to study the fixed point property \cite{LM2}.

Notice that, in the finite case, every time we apply a non-trivial retraction, we obtain a finite poset with fewer elements. Then, by induction, we can always obtain a minimal retract independently of the family of retractions that we are considering.

When passing to the infinite case this is no longer true. A rigorous definition of what constitutes a dismantling in the general case can be read in \cite{Kuk, Sch}. In these cases, we not only require a sequence (maybe indexed by an ordinal arbitrarily large) of composable retractions, but also some kind of compatibility condition that ensures us that the process has an end. At least we ask that for every point $x\in P$ the sequence of its successive images is eventually constant, or from a dynamics perspective, every trajectory eventually reaches a fixed point. An example of this is Kukie\l a's result \cite[Theorem 3.7, Remark 4.2]{Kuk} that for a sequence $r_n:P_n\to P_{n+1}\subseteq P_n$ of \textit{comparative retractions} (i.e. $r_n(p)$ and $p$ are comparable in $P_n$ for every $p\in P_n$), if every trajectory stabilizes, then $P=P_0$ strong deformation retracts to the subposet $Q$ of stable points. If not every trajectory reaches a fixed point, the stable points do not give information about homotopy invariants of $P$, moreover $Q$ could be empty. In the first case $Q$ is the colimit of the $P_n$. Comparative retractions are homotopy equivalences. Indeed, if $r:P\to P'$ is a comparative retraction, then $ir\le f\ge 1_P$, where $i:P'\to P$ denotes the inclusion, and $f:P\to P$ maps a point $p$ to $\max\{p, ir(p)\}$. Since comparative retractions are homotopy equivalences, Theorem \ref{main} says that $\text{colim }P_n$ has the weak homotopy type of $P$, so our result can be seen as an alternative to study the case of non-stabilizing trajectories.

Comparative retractions turn out to be a good source of applications of Theorem \ref{main}. In this case the colimit has a more geometric interpretation. The trajectory $p, r_0(p), r_1r_0(p), \ldots$ of a point $p\in P$ can be seen as analogous to the gradient flow lines determined by a Morse function in a manifold. Indeed, a Morse function encodes deformations, and so does a comparative retraction of A-spaces. The trajectory of a point $p\in P$ lives in the original space $P$, and two points in $P$ represent the same class in the trajectory space $\text{colim }P_n$ if their trajectories intersect. 

\begin{ej} \label{fence}
    The one-sided infinite fence is the space $X$ with underlying set $\N$ and relations $1<2>3<\ldots$. If $X_n$ denotes the subspace of $X$ given by those points $m\ge n$, then $n$ is a beat point of $X_n$, so there is a comparative retraction $r_n:X_n\to X_{n+1}$. In this case the space of trajectories has a unique point, so $X$ is homotopically trivial. In fact, it is not hard to prove that $X$ is contractible by exhibiting a concrete null-homotopy for the identity map.
\end{ej}

\begin{ej}
Consider the countable locally finite $T_0$ space $X$ on the left of Figure \ref{grilla}. There is a sequence $(r_n)_{n\in \N}$ of comparative retractions $r_n:X_n \to X_{n+1}$, with $X_1=X$. Each retraction removes a unique beat point. The order the points are removed, and the point each beat point is mapped to, are indicated by numbers and by red arrows in the right figure. 
\begin{figure}[ht] 
\begin{center}
\includegraphics[scale=0.35]{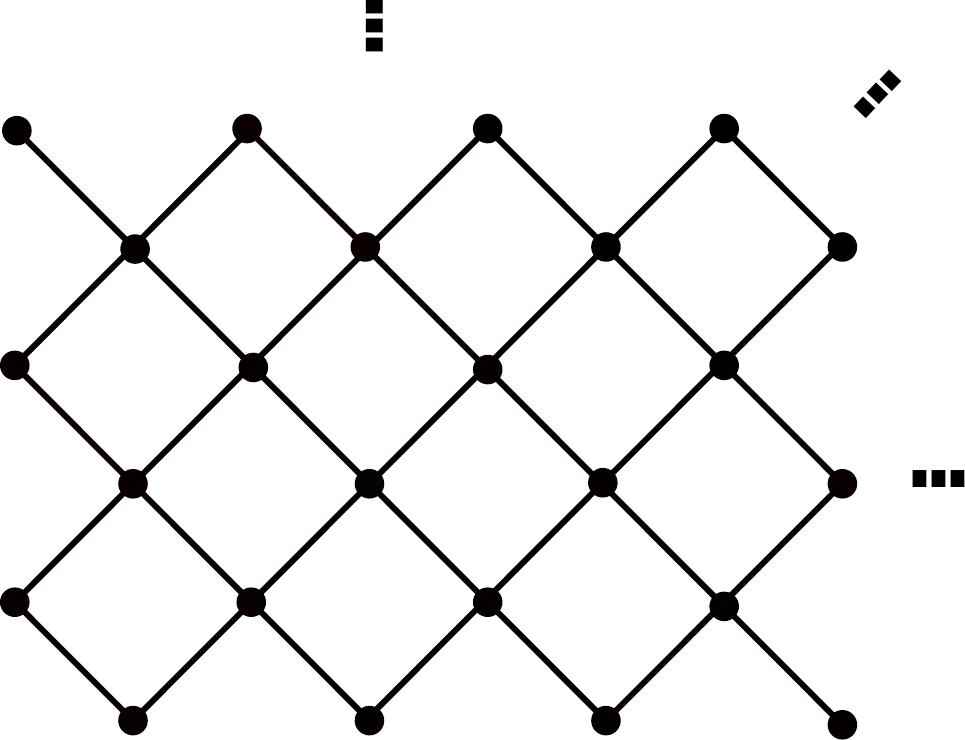} \qquad  \qquad \includegraphics[scale=0.35]{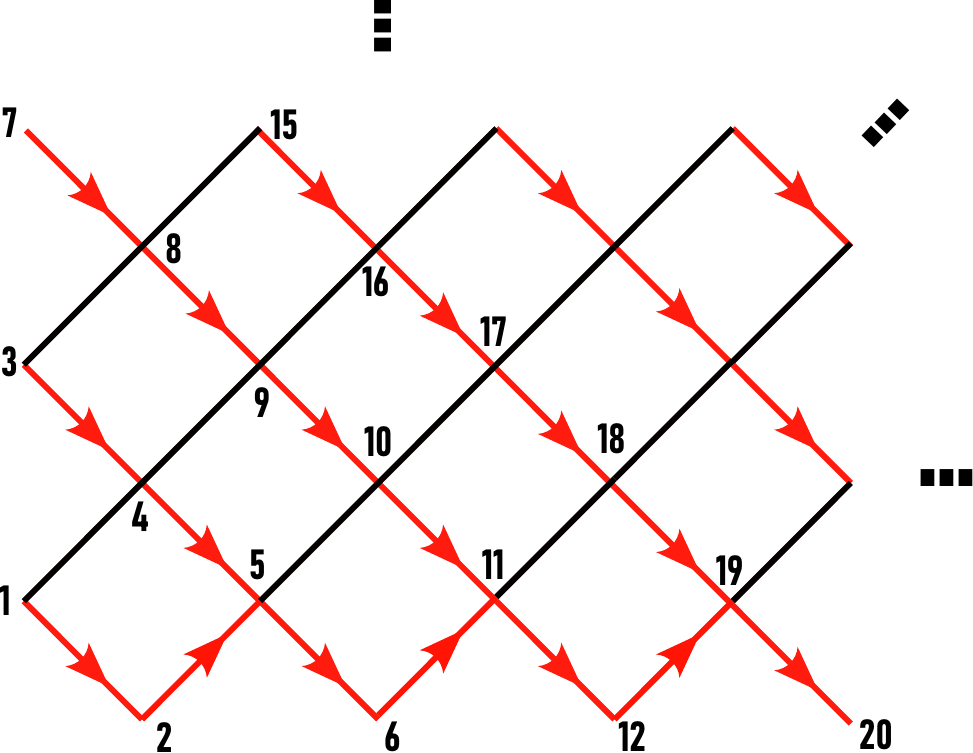}
\caption{The space $X$ on the left and a dismantling on the right.}\label{grilla}
\end{center}
\end{figure}
Since every pair of trajectories intersect, by Theorem \ref{main}, $X$ is homotopically trivial. This can also be proved by other methods, and, in fact, $X$ is contractible.
\end{ej}

One could ask if in the case of comparative retractions $\text{colim }P_n$ is homotopy equivalent to $P$. Example \ref{ejl} shows that this is not the case.






\begin{ej} \label{ejl}
The A-space that we will consider is named $L$ because it is a strong deformation retract of the Li Tower \cite{Li}.

The underlying set of $L$ is $$\{x_n\}_{n\in\N}\sqcup \{y\}\sqcup \{z_n\}_{n\in\N}$$. The topology is given by the following partial order, depicted in the left of Figure \ref{figlop}:
\begin{itemize}
	\item $x_i<x_j$ for every $i<j$.
	\item $z_i<x_j$ for every $i\le j$.
	\item $z_i<y$ for every $i$.
\end{itemize}

\begin{figure}[ht] 
\begin{center}
\includegraphics[scale=0.4]{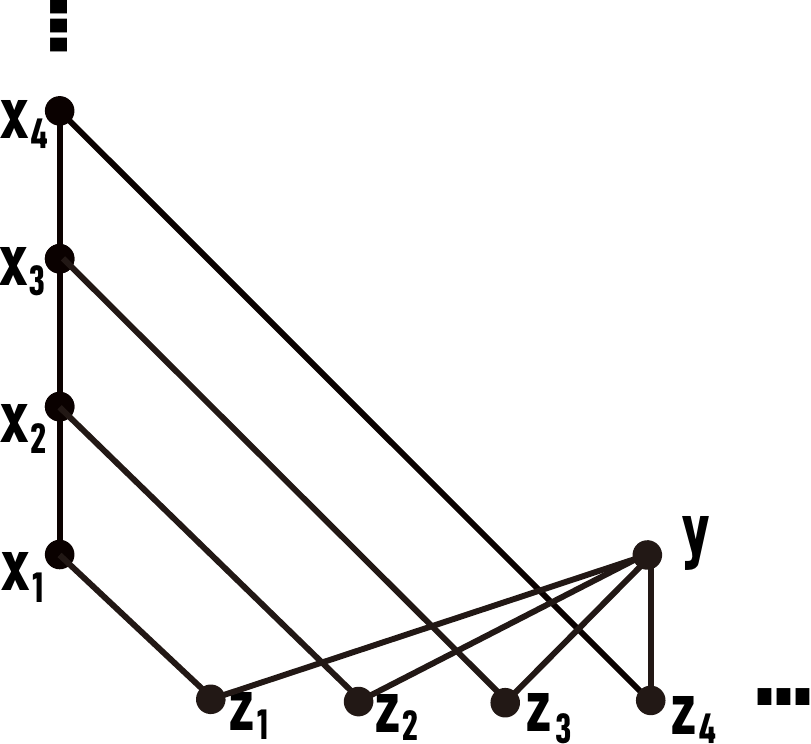} \qquad  \qquad \includegraphics[scale=0.4]{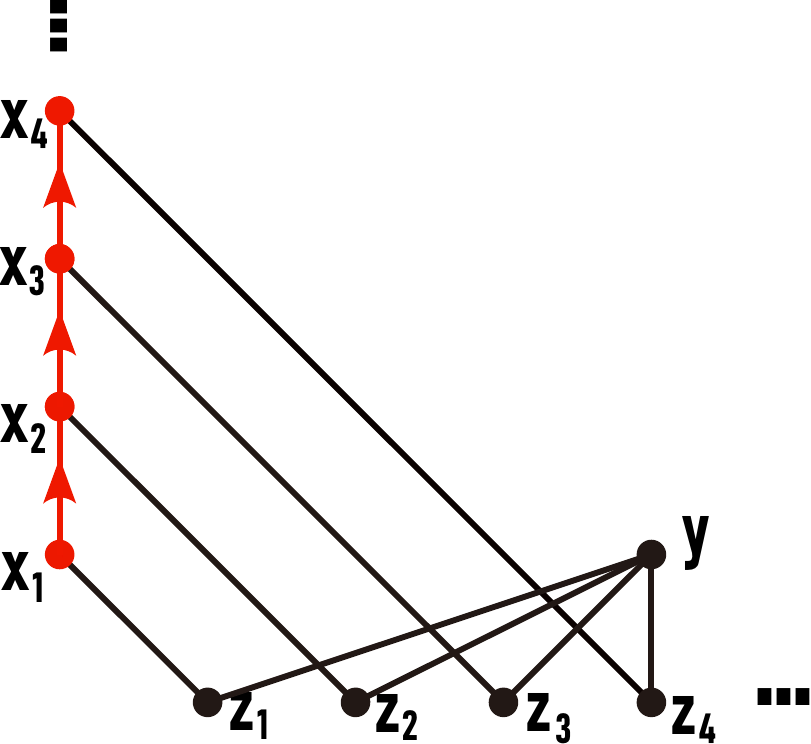}
\caption{The space $L$ and a dismantling.} \label{figlop}
\end{center}
\end{figure}


Define $L_1=L$ and recursively $L_{i+1}=L_i\smallsetminus \{x_i\}$ for $i\ge 1$. Since $x_i$ is an up-beat point of $L_i$, we have a comparative retraction $r_i:L_i\to L_{i+1}$ which maps $x_i$ to $x_{i+1}$. See the right of Figure \ref{figlop}.

If we consider the diagram $F:\N\to \tatop$ given by $F(n)=L_n$ and $F(n\to n+1)=r_n$, the colimit is the space of trajectories $L/R$ where two points are identified if their trajectories intersect. In this case there is only one orbit which is not a fixed point: $\{x_i\}_{i\in \N}$.
The topology of the trajectory space $\text{colim }F$ can be understood through Lemma \ref{diez}: since $L_i$ is a subspace of $L$ for every $i\in \N$, two classes $u,v$ in $\text{colim }F$ satisfy $u\le v$ if and only if there are representatives $\widetilde{u}, \widetilde{v}$ in $L$ such that $\widetilde{u}\le \widetilde{v}$. Using this it is easy to see that $\text{colim }F$ is homeomorphic to the poset in Figure \ref{colimitlop}. Once again this is the non-Hausdorff suspension of a countable discrete space.

\begin{figure}[ht] 
\begin{center}
\includegraphics[scale=0.4]{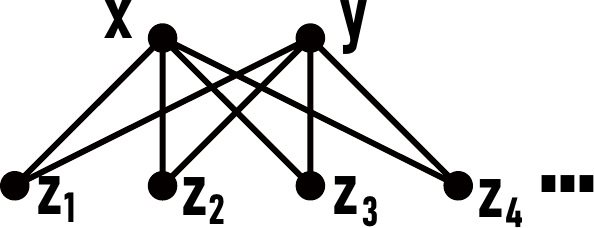}
\caption{The space of trajectories $L/R=\text{colim }F$.}\label{colimitlop}
\end{center}
\end{figure}


The associated polyhedron $|\mathcal{K}(\text{colim }F)|$ is a wedge of countably many copies of $S^1$ and, by Theorem \ref{main}, we have $\pi_1(L)\simeq \pi_1(\text{colim }F)\simeq \pi_1(\underset{n\in\N}{\bigvee}S^1)\simeq \underset{n\in\N}{*}\Z$.

Using this it is easy to check that $\text{colim }F$ is a weak core of $L$. Indeed, if $\text{colim F}\xrightarrow{q}Z$ is a quotient which identifies two distinct points $a_1,a_2\in \text{colim }F$, then there is a non-trivial loop $\gamma$, represented by a four element crown $b_1<b_2>b_3<b_4>b_1$ which contains both $a_1$ and $a_2$. Then $[q \gamma]=0 \in \pi_1(Z)$. This implies that $q_*:\pi_1(\text{colim }F)\to \pi_1(Z)$ is not an isomorphism, so $q$ is not a weak equivalence.


We claim that $p_1: L\to \text{colim }F$ is not a homotopy equivalence. Moreover, $\text{colim }F$ and $L$ do not have the same homotopy type. If $f:\text{colim }F \to L$ is a continuous map, then $f$ factorizes through the down-set (open set) $U_{f(x)}\cup U_{f(y)}$ generated by $f(x)$ and $f(y)$. But for any points $a,b\in L$, $U_a\cup U_b\subseteq L$ deformation retracts to a finite space, possibly by removing beat points $z_i$ for $i$ large enough. Since $f_*:\pi_1(\text{colim }F)\to \pi_1(L)$ factorizes through $\pi_1(U_{f(x)}\cup U_{f(y)})$, which is finitely generated, it cannot be an isomorphism. Thus, $f$ is not a homotopy equivalence.
\end{ej}

\section{A space with many weak cores} \label{sectionmany}

In this section we give an example of a countable A-space $X$ which has uncountably many (pairwise non-isomorphic) weak cores. The underlying set of $X$ is $\Z\cup \mathcal{P}_2(\Z)$, where $\mathcal{P}_2(\Z)$ denotes the set of 2-element subsets of $\Z$. The elements of $\Z$ are minimal in $X$, the points in $\mathcal{P}_2(\Z)$ are maximal, and $a\in \Z$ is smaller than $S\in \mathcal{P}(\Z)$ if $a\in S$. Note that $\mathcal{K}(X)$ is a graph and has the homotopy type of a wedge of countably many copies of $S^1$. 

\medskip

Let $Y$ be an A-space with underlying set a disjoint union of two copies $\Z$ and $\Z'=\{n'\vert \ n\in \Z\}$ of $\Z$. All the points in $\Z$ are minimal, all the points in $\Z'$ are maximal. Suppose that for any $n'\in \Z'$ the set $\hat{U}_{n'}=\{y\in Y \vert \ y<n'\}$ is an infinite subset of $\Z$. Moreover, $\hat{U}_{n'}=\Z$ for every $n< 0$. We prove now that under these conditions there exists a quotient map $q:X\to Y$ which is a weak equivalence.

The map $q:X\to Y$ we will construct is the identity $\Z \to \Z$ in the minimal sets. For each $n'\in \Z'$, $q^{-1}(n')$ will be a subset of $\mathcal{P}_2(\Z)\subseteq X$. We will define the sets $q^{-1}(n')$ for each $n'$ following this order: $0',1',-1',2',-2',\ldots$.

We first choose an ordering $S_1,S_2,\ldots$ of the elements of $\mathcal{P}_2(\Z)$. 

By assumption $\hat{U}_{0'}\subseteq \Z$ is infinite. Choose a numbering $a_1,a_2, \ldots$ of the elements in $\hat{U}_{0'}$. Define $q^{-1}(0')=\{\{a_i,a_{i+1} \}\vert  \ i\in \N \}\}$. Note that $q^{-1}(U_{n'})$ is a one-sided infinite fence $a_1<\{a_1,a_2\}>a_2<\{a_2,a_3\}>\ldots$, which is homotopically trivial by Example \ref{fence}.

Now we want to define $q^{-1}(1')$ in a similar way. We need to be careful since $q^{-1}(1')$ cannot intersect $q^{-1}(0')$. We choose an ordering $b_1,b_2,\ldots$ of $\hat{U}_{1'}$ in such a way that $\{b_i,b_{i+1}\}\neq \{a_j,a_{j+1}\}$ for any $i,j\in \N$. To achieve this, we can begin with any preliminary ordering $\mathfrak{b}_1,\mathfrak{b}_2,\ldots$ of $\hat{U}_{1'}$ and then reorder the elements as follows. Set $b_1=\mathfrak{b}_1$. Since $\mathfrak{b}_1$ and $\mathfrak{b}_2$ appear in only finitely many pairs $\{a_i,a_{i+1}\}$ and $\hat{U}_{1'}$ is infinite, there exists $k_1\neq 1,2$ such that $\{\mathfrak{b}_1,\mathfrak{b}_{k_1}\}$ and $\{\mathfrak{b}_{k_1},\mathfrak{b}_2\}$ are different from any $\{a_i,a_{i+1}\}$. Set $b_2=\mathfrak{b}_{k_1}$ and $b_3=\mathfrak{b}_{2}$. Let $l_2$ be the minimum of $\N\smallsetminus \{1,k_1,2\}$. We do with $\mathfrak{b}_2$ and $\mathfrak{b}_{l_2}$ the same we did with $\mathfrak{b}_1$ and $\mathfrak{b}_2$. Namely, since $\mathfrak{b}_2$ and $\mathfrak{b}_{l_2}$ appear only in finitely many $\{a_i,a_{i+1}\}$, we can choose $k_2\neq 1,k_1,2,l_2$ such that $\{\mathfrak{b}_2,\mathfrak{b}_{k_2}\}$ and $\{\mathfrak{b}_{k_2},\mathfrak{b}_{l_2}\}$ are different from any $\{a_i,a_{i+1}\}$. Set $b_4=\mathfrak{b}_{k_2}$ and $b_5=\mathfrak{b}_{l_2}$. Let $l_3$ be the minimum of $\N \smallsetminus \{1,k_1,2,k_2,l_2\}$. We continue in this way to obtain the desired ordering $b_1,b_2,\ldots$ of $\hat{U}_{1'}$. Once again, $q^{-1}(U_{1'})$ is a one-sided infinite fence.

If we continue in this way, defining $q^{-1}(n')$ for each $n\in \Z$ just being careful that each of these does not intersect the previously defined sets, we could have the problem that sets $q^{-1}(n')$ with $n\in \Z$ do not cover $\mathcal{P}_2(\Z)$. This is were we will use the ordering $S_1,S_2,\ldots$ of $\mathcal{P}_2(\Z)$ and the fact that $\hat{U}_{n'}=\Z$ for every $n<0$.

We have already defined $q^{-1}(0')$ and $q^{-1}(1')$. Let $m$ be the minimum $i\in \N$ such that $S_i$ is different from $\{a_j,a_{j+1}\}$ and from $\{b_j,b_{j+1}\}$ for any $j$. We now define an ordering $c_1,c_2,\ldots$ of $\hat{U}_{-1'}=\Z$ in such a way that 1. $\{c_1,c_2\}=S_m$ and 2. $\{c_i,c_{i+1}\}$ is different from $\{a_j,a_{j+1}\}$ and $\{b_j,b_{j+1}\}$ for any $i,j$.

We then proceed with $2'$, ignoring the order $S_1,S_2,\ldots$ and then we continue with $-2'$, defining an ordering of $\hat{U}_{-2'}=\Z$ in such a way that the first two elements in this ordering coincide with $S_{m'}$, the first pair in the ordering of the $S_i$ which has not appeared before.

In this way the sets $q^{-1}(n')$ with $n\in \Z$ form a partition of $\mathcal{P}_2(\Z)$, so we have a well defined function $q:X\to Y$.

It is easy to see that $q$ is continuous: If $a<\{a,b\}\in X$, then $q(\{a,b\})=n'$ for some $n\in \Z$ and $q(a)=a\in \hat{U}_{n'}$ by construction of $q$.

Also, $q$ is a quotient map. Indeed, on the one hand it is surjective. On the other hand, if $A\subseteq Y$ is not open, there exists $a<n' \in Y$ such that $a\notin A \ni n'$. Then, by construction of $q$, there exists $a\neq b<n'$ such that $\{a,b\}\in q^{-1}(n') \subseteq q^{-1}(A)$. Since $a\notin q^{-1}(A)$, $q^{-1}(A)$ is not open. 

Finally, $q:X\to Y$ is a weak equivalence by McCord's Theorem \ref{mccord} as $q^{-1}(U_y)$ is homotopically trivial for every $y\in Y$.

\medskip

In order to prove that $X$ has uncountably many weak cores, we need a condition that ensures that a space $Y$ as before is a weak core.

Let $Y$ be a space as above, but moreover assume that for each $n\in \Z$, $\Z \smallsetminus \hat{U}_{n'}$ is finite. Suppose $p:Y\to Z$ is a quotient map. We will prove that if $p$ is not injective, it cannot be a weak equivalence.

Suppose $p(n')=p(m')$ for some $n\neq m\in \Z$. Since both $\Z \smallsetminus \hat{U}_{n'}$ and $\Z \smallsetminus \hat{U}_{m'}$ are finite, there exist $a\neq b\in \Z$ such that both $a$ and $b$ are smaller than both $n'$ and $m'$. The crown $a<n'>b<m'>a$ represents a non-trivial element of $\pi_1(Y)$. Indeed, it induces a cycle $\gamma$ in the graph $|\mathcal{K}(Y)|$, so $\mu \gamma$ represents a non-trivial element of $\pi_1(Y)$. Since $p(n')=p(m')$, $p(\{a,n',b,m'\})$ has a maximum, so it is contractible. Thus $p_*([\mu \gamma])=0\in \pi_1(Z)$. This proves that $p$ is not a weak equivalence.

Suppose instead that $p(a)=p(b)$ for some $a\neq b \in \Z$. Then $a<-1'>b<-2'>a$ represents a non-trivial cycle, while $p(\{a,-1',b,-2'\})$ has a minimum, and thus is contractible. Again, $p$ is not a weak equivalence in this case either. 


Finally, if $p(a)=p(n')$ for some $a,n\in \Z$, take two different elements $b,c$ in $\hat{U}_{n'}$. Then $a<-1'>b<n'>c<-2'>a$ is a crown of six points which represents a non-trivial element of $\pi_1(Y)$. Every point in $p(\{a,-1',b,n',c,-2'\})$ is comparable with $p(a)=p(n')$, thus it is contractible, and $p$ is not a weak equivalence.

\medskip

To finish our example we exhibit uncountably many non-isomorphic spaces $Y$ which satisfy the conditions above.

Let $(u_n)_{n\in \N}\in \{0,1\}^{\N}$ be any sequence. We define the space $Y$ with underlying set $\Z \cup \Z'$, where all the points in $\Z$ are minimal and the remaining are maximal. As required we set $\hat{U}_{n'}=\Z$ for $n<0$. Also, set $\hat{U}_{0'}=\Z \smallsetminus \{0\}$. For $n\ge 1$ define $\hat{U}_{n'}=\Z \smallsetminus \{n-1,n\}$ if $u_n=0$ and $\hat{U}_{n'}=\Z \smallsetminus \{-1, n-1,n\}$ if $u_n=1$. Then $Y$ satisfies the required condition that $\Z\smallsetminus \hat{U}_{n'}$ is finite for every $n\in \Z$. 

Moreover, note that if $(u_n)_{n\in \N}, (v_n)_{n\in \N}\in \{0,1\}^{\N}$ are different, then the corresponding spaces $Y,Z$ are not homeomorphic. A homeomorphism $Y\to Z$ would induce an isomorphism between the bipartite complements of the bipartite graphs (comparability graphs) associated to $Y$ and $Z$, but these are non-isomorphic if the sequences are different. Since there are uncountably many sequences of $0$'s and $1$'s, $X$ has uncountably many weak cores. 

\section{Subdivisions and common approximations} \label{subdiv}

Our aim in this section is to prove that given A-spaces $X$ and $Y$ of the same weak homotopy type, there exists another A-space $Z$ and quotient maps $Z\xrightarrow{r}X$ and $Z\xrightarrow{r'} Y$ that are also weak equivalences. For this we introduce a notion of subdivision of A-spaces and prove that there is a quotient map from the subdivision to the original space, which is also a weak equivalence. This generalizes the result for barycentric subdivisions by Hardie and Vermeulen \cite[Proposition 2.2]{HV}.

\begin{defi}
	Let $X$ be a $T_0$ A-space and $K=\mathcal{K}(X)$ its associated simplicial complex. Let $L$ be a subdivision of $K$. That is, the vertex set of $L$ lies in $|K|$, for any simplex $\tau \in L$ there is a closed simplex $|\sigma|$ of $|K|$ which contains $\tau$, and the linear map $|L|\to |K|$ is a homeomorphism. We define the \textit{subdivision of $X$ associated to $L$} by $\mathcal{X}(L)$. We also define a map $\zeta :\mathcal{X}(L)\to X$ as follows. If $\tau \in \mathcal{X}(L)$, there exists a minimum simplex $\sigma_{\tau} \in \mathcal{K}(X)$ such that $\tau \subseteq |\sigma_{\tau}|$, which is called the carrier of $\tau$. Define $\zeta (\tau)=\max \sigma_{\tau}$.
\end{defi}



\begin{lema} \label{zeta}
	The map $\x (L)\xrightarrow{\zeta} X$ is a quotient map and a weak homotopy equivalence.
\end{lema}
\begin{proof}
	It is clear that $\zeta$ is order-preserving and surjective. To see that $\zeta$ is a quotient map take $A\subseteq X$ such that $\zeta^{-1}(A)\subseteq \mathcal{X}(L)$ is open. Let $x\in A$ and let $y<x$. Then $\{y,x\}$ is a simplex in $\mathcal{K}(X)$. The subdivision $L\restr_{\{y,x\}}$ of $\{y,x\}$ induced by $L$ is a subcomplex of $L$ with vertices $v_0=y, v_1, \ldots, v_k=x$ and simplices $\{v_i\}, \{v_i,v_{i+1}\}$. Since $\{y,x\}$ is the carrier of $\{v_0,v_1\}$, $\zeta (\{v_0,v_1\})=\max \{y,x\}=x\in A$. Since $\zeta^{-1}(A)$ is open and $\{v_0\}\le \{v_0,v_1\}$, $\{v_0\}\in \zeta^{-1}(A)$, so $y=\zeta (\{v_0\})\in A$. This proves that $A$ is open.
    

    By the McCord Theorem \ref{mccord}, to prove that $\zeta$ is a weak equivalence, it suffices to show that for every $x\in X$, $\zeta^{-1}(U_x)$ is homotopically trivial.

    A point $\tau\in \x(L)$ lies in $\zeta^{-1}(U_x)$ if and only $\max \sigma_{\tau}\le x$, or, equivalently, if $\sigma_{\tau} \subseteq U_x$, which amounts to $\sigma_{\tau} \in \kp (U_x)$. This, in turn, is equivalent to $\tau$ being in the subdivision $L\restr _{\kp(U_x)}$ of $\kp (U_x)$ induced by $L$. Thus $\zeta^{-1}(U_x)=\x(L\restr _{\kp(U_x)})$. Since $\kp (U_x)$ is a cone with apex $x$, $|\kp(U_x)|$ is contractible, and therefore so is $|L\restr _{\kp(U_x)}|$. This proves that $\zeta^{-1}(U_x)$ is homotopically trivial. 	
\end{proof}    

An alternative way to show that $\zeta$ is a weak equivalence is to notice that $\kp (\zeta)$ is an approximation to the identity.
    

	
	

\begin{teo}
	Let $X$ and $Y$ be A-spaces of the same weak homotopy type, then there exists an A-space $Z$ with quotient maps $Z\to X$ and $Z\to Y$ that are also weak homotopy equivalences.
\end{teo}
\begin{proof}
By adding topologically indistinguishable points to the common approximation $Z$ of $Q(X)$ and $Q(Y)$, we can obtain an approximation to $X$ and $Y$. Thus, we may assume $X$ and $Y$ are $T_0$. Let $f:|\mathcal{K}(X)|\to |\mathcal{K}(Y)|$ be a homotopy equivalence. By the simplicial approximation theorem, there exists a subdivision $L$ of $\mathcal{K}(X)$ and a simplicial map $\varphi: L\to \mathcal{K}(Y)$ such that $|\varphi|\simeq f$. We consider Spanier's simplicial mapping cylinder $M_{\varphi}$ \cite[Chapter 3, Exercise E]{Spa}. Choose a total ordering for the vertex set of $L$. The simplices of $M_{\varphi}$ are the subsets of sets $\{v_0,v_1, \ldots, v_k, \varphi(v_k), \varphi(v_{k+1}), \ldots , \varphi(v_l)\}$ where $v_0<v_1<\ldots <v_l$ and $\{v_0,v_1,\ldots, v_l\}$ is a simplex of $L$, and also the simplices of $\kp (Y)$. The realization of this complex is not homeomorphic to the topological mapping cylinder of $|\varphi|$, in general \cite{Coh}. However, we have canonical inclusions $i:L\to M_{\varphi}$ and $j:\kp (Y)\to M_{\varphi}$ and a simplicial retraction $r:M_{\varphi}\to \kp (Y)$ which maps $v\in L$ to $\varphi(v)$. Then $ri=\varphi$ and $1_{|M_{\varphi}|}\simeq |j||r| \text{ rel } |\kp (Y)|$. There is a homotopy $H: 1_{|M_{\varphi}|}\simeq |j||r| \text{ rel } |\kp (Y)|$ defined as follows. If $\alpha =\sum\limits_{i=0}^k t_i v_i +\sum\limits_{i=k}^l s_i \varphi(v_i) \in|M_{\varphi}|$, then $$H(\alpha,t)=t_0v_0+\ldots+t_{i-1}v_{i-1}+(t_i+t_{i+1}+\ldots+t_k-t)v_i+(t-t_{i+1}-t_{i+2}-\ldots-t_k)\varphi(v_i)+$$ $$+t_{i+1}\varphi(v_{i+1})+\ldots + t_{k-1}\varphi(v_{k-1})+(t_k+s_k)\varphi(v_k)+s_{k+1}\varphi(v_{k+1})+\ldots+s_l\varphi(v_l)$$ if $t_{i+1}+t_{i+2}+\ldots+t_k\le t\le t_{i}+t_{i+1}+\ldots+t_k$, and $H(\alpha,t)=|r|(\alpha)$ if $t\ge t_0+t_1+\ldots+t_k$. 

Since $|\varphi|$ and $|r|$ are homotopy equivalences, so is $|i|$. Then there exists a continuous retraction $p:|M_{\varphi}|\to |L|$ which is a homotopy inverse of $|i|$ (CW-pairs have the homotopy extension property).

Once again, by the simplicial approximation theorem there exists a subdivision $M$ of $M_{\varphi}$ and simplicial approximations $\pi: M\to L$ and $\rho:M\to \kp (Y)$ to $p$ and $|r|$.


\[\begin{tikzcd}
	& {|M|} & \\
	& {|M_\varphi|} \\
	\\
	{|L|} && {|\mathcal{K}(Y)|}
	\arrow[from=1-2, to=2-2]
	\arrow["{|\pi|}", curve={height=24pt}, from=1-2, to=4-1]
	\arrow["{|\rho|}"', curve={height=-24pt}, from=1-2, to=4-3]
	\arrow["p", shift left=3, from=2-2, to=4-1]
	\arrow["{|r|}"', from=2-2, to=4-3]
	\arrow["{|i|}", from=4-1, to=2-2]
	\arrow["{|\varphi|}", from=4-1, to=4-3]
\end{tikzcd}\]

We claim that $\pi$ is surjective on simplices. Indeed, if $\tau$ is a simplex of $L$, then the subdivision $M\restr _\tau$ of $\tau$ induced by $M$ is a triangulation of $\tau$, and the fact that $p$ is a retraction and $\pi$ approximates $p$, implies that $\pi: M\restr_\tau \to \tau$ is a Sperner coloring. By the Sperner lemma \cite{Spe} a simplex of $M\restr_\tau$ is mapped bijectively to $\tau$ by $\pi$. 

We prove now that $\x(\pi):\x(M)\to \x(L)$ is a quotient. We have already shown it is surjective. This automatically implies that it is a quotient. Namely, suppose $A\subseteq \x(L)$ and $\x(\pi)^{-1}(A)\subseteq \x(M)$ is open. Let $\tau \in A$ and let $\tau' \subseteq \tau$. Since $\pi$ is surjective on simplices, there exists $\sigma \in M$ such that $\pi (\sigma)=\tau$. Then there exists $\sigma' \subseteq \sigma$ such that $\pi (\sigma')=\tau'$. Since $\sigma \in \x(\pi)^{-1}(A)$, which is open, then $\sigma'\in \x(\pi)^{-1}(A)$ and $\tau'\in A$. This proves that $A$ is open, so $\x(\pi)$ is a quotient.

As $p$ is a homotopy equivalence, so is $|\pi|$. Thus $\x(\pi)$ is a weak equivalence.

Since $|r|$ is a retraction and a homotopy equivalence, the same argument shows that $\x(\rho)$ is a quotient and a weak equivalence.
	
Finally, Lemma \ref{zeta} implies that $\zeta \mathcal{X}(\pi):\x(M)\to X$ and $\zeta \mathcal{X}(\rho):\x(M)\to Y$ are quotient maps and weak equivalences. 
\end{proof}

As an immediate consequence we get the next result.
\begin{coro}
	If $X_0,X_1,\dots, X_n$ are A-spaces of the same weak homotopy type, there exists an A-space $Z$ with quotient maps $Z\xrightarrow{r_i}X_i$ that are also weak homotopy equivalences.
\end{coro} 

\begin{obs}
    The common approximation $Z$ of $X$ and $Y$ can be taken to be $T_0$, even if $X$ or $Y$ are not. Concretely, if $Z$ is an A-space, there exists a $T_0$ A-space $T(Z)$ and a quotient map $q:T(Z)\to Z$ which is a homotopy equivalence.

    Consider the equivalence relation on $Z$ defined by $z\sim z'$ if $z\le z'$ and $z'\le z$. From each class $c\in Z/_\sim=Q(Z)$ we choose an element $s(c)\in c$. Let $\hat{c}$ be the set obtained from $c$ by adding a new element $m_c$. Define a total order on $\hat{c}$ in such a way that $s(c)$ is the minimum and $m_c$ is the maximum. Let $T(Z)$ be the union of the sets $\hat{c}$. That is $T(Z)=Z\cup \{m_c\}_{c\in Q(Z)}$. We define an order $\trianglelefteq$ on $T(Z)$ by keeping the total order in each set $\hat{c}$ and setting $z' \trianglelefteq z$ for each $z\in \hat{c}, z'\in \hat{c}'$ if $c'< c \in Q(Z)$ (i.e. if $c'\neq c$ and any element of $c'$ is smaller than any element of $c$ in $Z$). Then $\trianglelefteq$ is indeed a partial order, so $T(Z)$ is $T_0$. Define $q:T(Z)\to Z$ by $q(z)=z$ if $z\in Z$ and $q(m_c)=s(c)$. Then $q$ is order-preserving. Let $f:Z\to T(Z)$ be the map that maps $z\in Z$ to $s([z])$. Then $f$ is also order-preserving. Since $qf\le 1_{Z}$ and $fq\le 1_{T(Z)}$, $q$ is a homotopy equivalence.  

    We prove $q$ is a quotient. It is clearly surjective. Let $A\subseteq Z$ be such that $q^{-1}(A)\subseteq T(Z)$ is open. Let $z\in A$ and let $z'\le z$. If $[z']\neq [z]$, $z'\trianglelefteq z \in q^{-1}(A)$, so $z'\in q^{-1}(A)$ and then $z'\in A$. If $[z']=[z]$, then $s([z])\trianglelefteq z \in q^{-1}(A)$, so $s([z])\in q^{-1}(A)$. Thus $m_{[z]}\in q^{-1}(A)$ and since $z'\trianglelefteq m_{[z]}$, $z'\in q^{-1}(A)$. Thus $z'\in A$. This proves that $A$ is open. 
\end{obs}

\end{document}